\documentclass{amsart}
\usepackage{amsmath, amsxtra, amsfonts, amssymb, amsthm, comment, enumerate, mathtools}
\usepackage[final]{hyperref}

\allowdisplaybreaks

\newcommand{\R}[0]{\mathbb R}
\newcommand{\ml}[0]{\mathcal L}
\newcommand{\vnorm}[1]{\left \lvert #1 \right \rvert}
\newcommand{\Vnorm}[1]{\left \lVert #1 \right \rVert}

\theoremstyle{theorem}
\newtheorem{thm}{Theorem}
\newtheorem{prp}{Proposition}

\DeclareMathOperator{\diverg}{div}

\makeatletter
\renewcommand\subsection{\@startsection{subsection}{2}
  \z@{.5\linespacing\@plus.7\linespacing}{-.5em}
  {\normalfont\itshape}}
\makeatother

\begin{document}

\author{Lorenzo Riva \and Nathan Pennington}

\address{Creighton University}
\email[Lorenzo Riva]{LorenzoRiva@creighton.edu}
\email[Nathan Pennington]{NathanPennington@creighton.edu}

\title{Low regularity of non-$L^2(R^n)$ local solutions to gMHD-alpha systems}
\date{05/26/2020}

\maketitle

\begin{abstract}
	The Magneto-Hydrodynamic (MHD) system of equations governs viscous fluids subject to a magnetic field and is derived via a coupling of the Navier-Stokes equations and Maxwell's equations. It has recently become common to study generalizations of fluids-based differential equations. Here we consider the generalized Magneto-Hydrodynamic alpha (gMHD-$\alpha$) system, which differs from the original MHD system by the presence of additional non-linear terms (indexed by the choice of $\alpha$) and replacing the Laplace operators in the equations by more general Fourier multipliers with symbols of the form $-\vnorm{\xi}^\gamma / g(\vnorm{\xi})$. In \cite{penn1}, one of the authors considered the problem with initial data in Sobolev spaces of the form $H^{s,2}(\mathbb{R}^n)$  with $n \geq 3$. Here we consider the problem with initial data in $H^{s,p}(\mathbb{R}^n)$ with $n \geq 3$ and $p > 2$, with the goal of minimizing the regularity required to obtain unique existence results.
\end{abstract}

\section{Introduction}

This paper is concerned with the generalized Magneto-Hydrodynamic alpha (gMHD-$\alpha$) system of equations, reported below in its full generality:
\begin{align}
	& \partial_t v + (u \cdot \nabla) v + \sum_{i = 1}^n v_i \nabla u_i - \nu_1 \ml_1 v + \frac{1}{2} \nabla \vnorm{B}^2 = - \nabla p + (B \cdot \nabla)B, \label{u-eq} \\
	& \partial_t B + (u \cdot \nabla) B - (B \cdot \nabla) u - \nu_2 \ml_2 B = 0, \label{b-eq} \\
	& v = (1 - \alpha^2 \ml_3) u, \label{v-eq} \\
	& \diverg u = \diverg B = 0, \label{div-eq} \\
	& u(0,x) = u_0(x), \quad B(0,x) = B_0(x), \quad x \in \mathbb R^n.	\label{init}
\end{align}
Since these equations govern the motion of fluids subject to a magnetic field, (almost) each term has a specific physical meaning: in the above, $u$ is the fluid velocity, $B$ the magnetic field, and $p$ the scalar-valued pressure of the fluid; for the constants, $\nu_1 > 0$ is the fluid viscosity, $\nu_2 > 0$ the magnetic diffusion, and $\alpha > 0$ a constant coming from varying the Hamiltonian that originally gave rise to the standard MHD equations (see \cite{linshiz}). Finally, the $\ml_i$ terms are Fourier multipliers with symbol $-\vnorm{\xi}^{\gamma_i} / g_i(\vnorm{\xi})$, where $g_i$ is a positive scalar function and $\gamma_i > 0$.

The standard MHD system is the special case obtained when setting $\alpha = 0$, $g_1 = g_2 = 1$, and $\gamma_1 = \gamma_2 = 2$, so that $v = u$ and $\ml_1 = \ml_2 = \Delta$. The existence of a global solution up to initial conditions is a classic result, at least in two dimensions. Unfortunately, the presence of nonlinear terms makes the MHD equations particularly complex to solve in arbitrary dimensions, and so a common strategy has been to study modified versions of them. 

One modification, termed Lagrangian Averaged MHD-$\alpha$ after the Lagrangian Averaged Navier-Stokes equation, is obtained from Equations (\ref{u-eq})-(\ref{init}) by setting $\ml_i = \Delta$ for $i = 1,2,3$ ($\gamma_i = 2$ and $g_i = 1$). Linshiz and Titi proved the existence of a global solution for smooth initial data in three dimension \cite{linshiz}. Another version is obtained by setting $\alpha = 0$ and $g_1 = g_2 = 1$ and leaving the $\gamma_i$'s unspecified. Zhao and Zhu used these generalized operators to guarantee a global solution to Equations (\ref{u-eq})-(\ref{init}) in the case of $n = 3$, provided that $g_1 = g_2 = g_3 = 1$, $\gamma_1 = \gamma_2 = n/2$, and $\gamma_3 = 2$ \cite{zhao}.

The first incorporation of a non-constant value for any $g_i$ appeared in \cite{tao}, where Tao proved the existence of a unique global solution to the generalized Navier-Stokes equation ($B_0 = \alpha = 0$) when $\gamma_1 = n/2 + 1$ and when $g_1$ is a radial non-decreasing function bounded below satisfying
\begin{equation} \label{tao-cond}
	\int_1^\infty \frac{ds}{s g_1(s)^4} = \infty,
\end{equation}
the prototypical example of which is essentially a logarithm.

Wu obtained a similar result for the generalized MHD system in \cite{wu}, specifically showing that there is a unique global solution provided that $u_0, B_0 \in H^{r,2}(\R^n)$ with $r > n/2 + 1$; $\gamma_1 \geq n/2 + 1$, $\gamma_2 > 0$, and $\gamma_1 + \gamma_2 \geq 1$; and $g_1, g_2$ are non-decreasing, bounded below by 1, and satisfy
\begin{equation} \label{wu-cond}
	\int_1^\infty \frac{ds}{s (g_1(s) + g_2(s))^2} = \infty.
\end{equation}
This work was ultimately extended to the gMHD-$\alpha$ system in \cite{yamazaki}, where Yamazaki obtained a unique global solution in three dimensions provided that $\gamma_1 + \gamma_2 + \gamma_3 \geq 5$, $\min \{\gamma_1, \gamma_3\} > \gamma_2 > 0$, $\gamma_3 + 2 \gamma_1 > 3$, and the $g_i$ satisfy
\begin{equation} \label{yama-cond}
	\int_1^\infty \frac{ds}{s g_1(s)^2 g_2(s) g_3(s)^2} = \infty.
\end{equation}

In \cite{penn1}, one of the authors considered a generalization of the equations in \cite{zhao} with the incorporation of non-constant $g_i$, $i = 1, 2, 3$, while still leaving $\ml_3 = \Delta$, and guaranteed a unique global solution. 

In this paper we will extend those results for the case of $\gamma_3 \neq 2$ and non-constant $g_3$. We will particularly focus on the case of low-regularity initial data in a non-$L^2(\R^n)$ setting to then obtain, in the future, global $L^p(\R^n)$ solutions using an interpolation technique, the details of which can be found in \cite{penn4} and \cite{gallagher}.

The rest of the paper is organized as follows. Section 2 is devoted to explaining the notation we will use and some supporting results necessary for the algorithm. Section 3 contains the main result (Theorem \ref{main-thm}) of this paper and its proof. We end this section with two important spacial cases of Theorem \ref{main-thm}.

\begin{thm} \label{special-case}
	Let $g_1, g_2, g_3 : [0, \infty) \to \R$ be non-decreasing functions bounded below by 1 satisfying
	\begin{equation} \label{mikhlin-condition}
		g_i^{(k)}(s) \leq Cs^{-k}
	\end{equation}
	for $i = 1,2,3$ and $0 \leq k \leq n/2+1$. Moreover, assume $0 \leq \gamma_3 \leq 1$ and $p, q \geq n$ with $2p > q$.

	Then, for any divergence-free $u_0 \in L^p(\R^n)$ and $B_0 \in L^q(\R^n)$, there exists a unique local solution $(u,B)$ to the generalized MHD-$\alpha$ system (\ref{u-eq})-(\ref{init}) provided that
	\begin{align*}
		\gamma_1^- & > 6 - \gamma_3, \\
		\gamma_2^- & > 1 + \frac{n}{p}.
	\end{align*}
\end{thm}

\noindent The condition (\ref{mikhlin-condition}) is a modification of the condition in the Mikhlin multiplier thoerem that is necessary for supporting estimates in Proposition (\ref{semigroup-est}) and (\ref{inverse-semigroup-est}). The functions that satisfy it are still essentially logarithms, the same type of functions that satisfy (\ref{tao-cond})-(\ref{yama-cond}).

\begin{thm} \label{special-case-2}
	Let $\gamma_3^- - 1 \leq \frac{n}{2p} \leq \gamma_3^-$, $\frac{n}{2q} - 1 + \gamma_3^- \leq \frac{n}{2p}$, and let $p, q \geq n$ with $q < 3p/2$. Moreover, assume that $g_1, g_2, g_3$ satisfy the inequality (\ref{mikhlin-condition}). Then, for each divergence-free $u_0 \in H^{n/2p,p}(\R^n)$ and $B_0 \in H^{n/2q,q}(\R^n)$, there exists a unique local solution $(u,B)$ to the generalized MHD system from Equations (\ref{u-eq})-(\ref{init}) provided that
	\begin{align*}
		\gamma_1^- & > 6 - \gamma_3^- - \frac{n}{p}, \\
		\gamma_2^- & > 1 + \frac{n}{2p}.
	\end{align*}
\end{thm}

Note that, in the statement of Theorem \ref{special-case} and \ref{special-case-2} and in what follows, we use $x^- = x - \varepsilon$ for some positive $\varepsilon$, i.e. $x^-$ denotes a number arbitrarily close to, but strictly smaller than, $x$.

\section{Notation and supporting facts}

We let $H^{r,p}(\R^n)$ be the usual Sobolev space, and we write $\Vnorm{f}_{r,p}$ to mean $\Vnorm{f}_{H^{r,p}(\R^n)}$ and $\Vnorm{f}_{p}$ for $\Vnorm{f}_{L^p(\R^n)}$. Due to the nature of the procedure we will use, we require that the solutions live in an auxiliary continuous-in-time space $C^T_{a;r,p}(\R^n)$ defined by
\begin{equation*}
	C^T_{a;r,p}(\R^n) := \left\{ f \in C((0,T), H^{r,p}(\R^n)) : \Vnorm{f}_{a;r,p} < \infty \right\},
\end{equation*}
where $T > 0$, $a \geq 0$, $C(X,Y)$ is the space of continuous maps $X \to Y$, and
\begin{equation*}
	\Vnorm{f}_{a;r,p} := \sup_{(0,T)} t^a \Vnorm{f(t)}_{r,p}.
\end{equation*}
Finally, we denote by $\dot{C}^T_{a;r,p}(\R^n)$ the subspace of $C^T_{a;r,p}(\R^n)$ consisting of functions $f$ such that \mbox{$\lim_{t \to 0^+} t^{a} f(t) = 0$} and by $BC(X,Y) \subset C(X,Y)$ the subspace of bounded continuous maps $X \to Y$.

The following are some supporting propositions that we will use throughout the paper. This first one is a product estimate, the proof of which can be found in Chapter 2 of \cite{taylor1}:
\begin{prp} \label{product-est}
	If $r \geq 0$ and $1 < p \leq \infty$, then
	\begin{equation*}
		\Vnorm{fg}_{r,p} \leq C \left( \Vnorm{f}_{p_1} \Vnorm{g}_{r,p_2} + \Vnorm{f}_{r, q_1} \Vnorm{g}_{q_2} \right),
	\end{equation*}
	where
	\begin{equation*}
		\frac{1}{p} = \frac{1}{p_1} + \frac{1}{p_2} = \frac{1}{q_1} + \frac{1}{q_2}
	\end{equation*}
	and $p_1, p_2, q_1, q_2 \in [1, \infty]$.
\end{prp}

The following is a useful Sobolev embedding which is a straightforward extension of a result from Chapter 13 of \cite{taylor2}:

\begin{prp} \label{sobolev-est}
	Let $s \geq r$ and $(s-r)p < n$. Then
	\begin{equation*}
		\Vnorm{f}_{r,q} \leq C \Vnorm{f}_{s,p}
	\end{equation*}
	provided that
	\begin{equation*}
		\frac{1}{q} - \frac{r}{n} = \frac{1}{p} - \frac{s}{n}.
	\end{equation*}
\end{prp}

Our next result follows from a simple calculus exercise:
\begin{prp} \label{integral-est}
	If $0 < a, b \in \R$, then
	\begin{equation*}
		\sup_{t \in [0,T]} \int_0^t (t-s)^{-a} s^{-b} \, ds \leq CT^{1 - a - b},
	\end{equation*}
	provided that $a + b < 1$.
\end{prp}

Our final two propositions consist of an estimate for the semigroup $e^{t\ml_i}$ analogous to similar results for the heat kernel $e^{t\Delta}$ and an estimate for the operator $(1 - \ml_i)^{-1}$. The proofs of both propositions can be found in \cite{penn3}. We recall that $x^-$ is a number arbitrarily close to, but strictly smaller than, $x$.
\begin{prp} \label{semigroup-est}
	Let $1 < p_1 \leq p_2 < \infty$, $r_1 \leq r_2$, $g(x)$ be a non-decreasing function bounded below by 1 satisfying $\vnorm{g^{(k)}(x)} \leq C \vnorm{x}^{-k}$ for all $1 \leq k \leq n/2 + 1$. Then
	\begin{equation*}
		e^{t\ml_i} : H^{r_1,p_1}(\R^n) \to H^{r_2, p_2}(\R^n)
	\end{equation*}
	and
	\begin{equation*}
		\Vnorm{e^{t\ml_i} f}_{r_2, p_2} \leq t^{-(r_2 - r_1 + n/p_1 - n/p_2)/\gamma_i^-} \Vnorm{f}_{r_1, p_1}. 
	\end{equation*}
\end{prp}

\noindent Note that it is this proposition which necessitates the requirements on the $g_i$'s.

\begin{prp} \label{inverse-semigroup-est}
	Let $1 < p < \infty$, $r \in \R$, $g(x)$ be a non-decreasing function bounded below by 1 satisfying $\vnorm{g^{(k)}(x)} \leq C \vnorm{x}^{-k}$ for all $1 \leq k \leq n/2 + 1$. Then
	\begin{equation*}
		\Vnorm{(1 - \ml_i)^{-1} f}_{r,p} \leq C \Vnorm{f}_{r - \gamma_i^-, p}.
	\end{equation*}
\end{prp}

\section{Main Theorem and Proof}

In this section we state the most general form of the theorem and then proceed with the proof.

\begin{thm} \label{main-thm}
	Let $g_1, g_2, g_3 : [0, \infty) \to \R$ be non-decreasing functions bounded below by 1 satisfying
	\begin{equation*}
		g_i^{(k)}(s) \leq Cs^{-k}
	\end{equation*}
	for $i = 1,2,3$ and $0 \leq k \leq n/2+1$. Let $r_0, r_1, r_2 \geq 0$ and let $p_0, p_1, p_2 \geq n$ with $p_0 \leq p_1$ and $p_2 < 2p_0$. Moreover, assume that
	\begin{gather*}
		\gamma_3^- - 1 \leq r_0 \leq \gamma_3^- \leq r_1, \\
	r_2 - 1 + \gamma_3^- \leq r_0, \\
	r_2 \leq r_0 < \frac{n}{p_1}, \\
	2r_1 \geq \max \left\{2, 1 + \gamma_3^- - \frac{n}{p_0} + \frac{2n}{p_1} \right\}, \\
	r_2 < \min \left\{\frac{n}{p_2}, \frac{2n}{p_2} - \frac{n}{p_0} \right\}.
	\end{gather*}

	Then, for any divergence-free $u_0 \in H^{r_0,p_0}(\R^n)$ and $B_0 \in H^{r_2,p_2}(\R^n)$, there exists a unique local solution $(u,B)$ to the generalized MHD-$\alpha$ system (\ref{u-eq})-(\ref{init}) provided that
	\begin{align*}
		\gamma_1^- & > 3r_1 - 2r_0 - \gamma_3^- + \frac{3n}{p_0} - \frac{3n}{p_1}, \\
		\gamma_1^- & > 1 - 2r_2 + r_1  - \gamma_3^- - \frac{n}{p_1} + \frac{2n}{p_2}, \\
		\gamma_2^- & > 1 - r_0 + \frac{n}{p_0}.
	\end{align*}
\end{thm}
\noindent Note: Theorem \ref{special-case} may be recovered by setting $r_0 = r_2 = 0$, $r_1 = 2$, $p := p_0 = p_1$, and $q := p_2$, while Theorem \ref{special-case-2} may be recovered by setting $r_0 = n/(2p)$, $r_1 = 2$, $r_2 = n/(2q)$, $p := p_0 = p_1$, and $q := p_2$.
\newline

\begin{proof}
For the sake of clarity and to highlight some technical details, the proof will be divided in subsections. We first write the generalized MHD-$\alpha$ system in a more helpful form. Without loss of generality, we set $\alpha = \nu_1 = \nu_2 = 1$. We pass to divergence-free vector fields by applying the Hodge operator $P$ to Equations (\ref{u-eq}) and (\ref{b-eq}) (more information about the Hodge operator can be found in Chapter 11 of \cite{rieusset}), and then we apply $(1 - \ml_3)^{-1}$ to Equation (\ref{u-eq}). By noting that $P$, $(1 - \ml_3)^{-1}$, and $\partial_t$ all commute since they are Fourier multipliers, we obtain
\begin{equation*}
	\partial_t u + P(1 - \ml_3)^{-1} \left((u \cdot \nabla) v + \sum_{i = 1}^n v_i \nabla u_i - (B \cdot \nabla) B \right) - \ml_1 u = P \left( -\nabla p - \frac{1}{2} \nabla \vnorm{B}^2 \right) = 0.
\end{equation*}
An application of the divergence-free condition in Equation (\ref{div-eq}) allows us to rewrite the terms of the form $(x \cdot \nabla) y$ as $\diverg(x \otimes y)$. Note that $x \otimes y$ is the matrix whose $(i,j)$ entry is $x_i y_j$, so that the product estimate in Proposition \ref{product-est} applies to $x \otimes y$. We then have the following system:
\begin{align*}
	& \partial_t u + P(1 - \ml_3)^{-1} \left(\diverg (u \otimes v) + \sum_{i = 1}^n v_i \nabla u_i - \diverg (B \otimes B) \right) - \ml_1 u = 0, \\
	& \partial_t B + P \left(\diverg (u \otimes B) - \diverg (B \otimes u)\right) - \ml_2 B = 0, \\
	& v = (1 - \ml_3) u, \\
	& \diverg u = \diverg B = 0, \\
	& u(x,0) = u_0(x), \quad B(0,x) = B_0(x), \quad x \in \mathbb R^n.
\end{align*}

An application of Duhamel's principle shows that $(u,B)$ is a solution to the system if and only if $(u,B)$ is a fixed point of the map $\Phi(u,B) := (\Phi_1(u,B), \Phi_2(u,B))$ defined by
\begin{align*}
	\Phi_1(u,B) & := e^{t\ml_1} u_0 - \int_0^t e^{(t-s) \ml_1} \left( W_1 (u,v) + W_2(u,v) - W_1 (B,B)\right) ds, \\
	\Phi_2(u,B) & := e^{t\ml_2} B_0 - \int_0^t e^{(t-s)\ml_2} \left( W_3 (u,B) - W_3(B,u) \right) ds,
\end{align*}
where
\begin{align*}
	W_1(x,y) & = P (1-\ml_3)^{-1} \diverg (x \otimes y), \\
	W_2(x,y) & = P (1-\ml_3)^{-1} \left( \sum_{i = 1}^n y_i \nabla x_i \right), \\
	W_3(x,y) & = P \diverg (x \otimes y).
\end{align*}

By the contraction mapping theorem, it suffices to show that $\Phi$ is a contraction on the space $X_{T,M} \times Y_{T,M}$, where
\begin{multline*}
	X_{T,M} := \bigg\{ f \in BC\left([0,T), H^{r_0,p_0}(\R^n)\right) \cap \dot{C}_{a_1, r_1, p_1}(\R^n) \\
	\text{ and } \sup_{(0,T)} \left( \Vnorm{f(t) - e^{t\ml_1}u_0}_{r_0, p_0} + \Vnorm{f(t)}_{a_1; r_1, p_1} \right) < M \bigg\}
\end{multline*}
and
\begin{equation*}
	Y_{T,M} := \bigg\{ f \in BC([0,T), H^{r_2,p_2}(\R^n)) \text{ and }
	\sup_{(0,T)} \Vnorm{f(t) - e^{t\ml_2}B_0}_{r_2, p_2} < M \bigg\}
\end{equation*}
for some $0 < T < 1$ and $M > 0$.

Following the methods in \cite{penn2,kato}, we will complete the proof by showing that
\begin{align*}
	I_1 & = \sup_{(0,T)} t^{a_1} \Vnorm{e^{t\ml_1} u_0}_{r_1,p_1} < M/4, \\
	I_2 & = \sup_{(0,T)} \Vnorm{\int_0^t e^{(t-s) \ml_1} \left( W_1 (u,v) + W_2(u,v) - W_1 (B,B)\right) ds}_{r_0,p_0} < M/4, \\
	I_3 & = \sup_{(0,T)} t^{a_1} \Vnorm{\int_0^t e^{(t-s) \ml_1} \left( W_1 (u,v) + W_2(u,v) - W_1 (B,B)\right) ds}_{r_1,p_1} < M/4, \\
	I_4 & = \sup_{(0,T)} \Vnorm{\int_0^t e^{(t-s) \ml_2} \left( W_3 (u,B) - W_3(B,u) \right) ds }_{r_2,p_2} < M/4.\\
\end{align*}

We start with $I_1$. If $\varphi$ is in the Schwartz space, we have that
\begin{align*}
	I_1 & = \sup_{(0,T)} t^{a_1} \Vnorm{e^{t\ml_1} (u_0 - \varphi + \varphi)}_{r_1,p_1} \\
	& \leq \sup_{(0,T)} t^{a_1} \Vnorm{e^{t\ml_1} (u_0 - \varphi)}_{r_1, p_1} + \sup_{(0,T)} t^{a_1} \Vnorm{e^{t\ml_1} \varphi}_{r_1, p_1} \\
	& \leq \sup_{(0,T)} t^{a_1} t^{-a_1} \Vnorm{u_0 - \varphi}_{r_0, p_0} + \sup_{(0,T)} t^{a_1} \Vnorm{\varphi}_{r_1, p_1} \\
	& \leq \Vnorm{u_0 - \varphi}_{r_0, p_0} + T^{a_1} \Vnorm{\varphi}_{r_1, p_1},
\end{align*}
provided that (by Proposition \ref{semigroup-est})
\begin{equation*}
	0 \leq a_1 = \frac{r_1 - r_0 + \frac{n}{p_0} - \frac{n}{p_1}}{\gamma_1^-} < 1
\end{equation*}
and $p_0 \leq p_1$. We can choose $\varphi$ so that $\Vnorm{u_0 - \varphi}_{r_0, p_0}$ is arbitrarily small, and then we can choose $T$ small enough to reduce $T^{a_1} \Vnorm{\varphi}_{r_1, p_1}$ so that the sum of the two is bounded by $M/4$.

\subsection{$I_2$ and $I_3$.}

Minkowski's inequality gives us
\begin{align*}
	I_2 & \leq J_1 + J_2 + J_3, \\
	I_3 & \leq K_1 + K_2 + K_3,
\end{align*}
where
\begin{align*}
	J_1 & := \sup_{(0,T)} \int_0^t \Vnorm{e^{(t-s) \ml_1} W_1 (u,v)}_{r_0,p_0} ds, & K_1 & := \sup_{(0,T)} t^{a_1} \int_0^t \Vnorm{e^{(t-s) \ml_1} W_1 (u,v)}_{r_1,p_1} ds, \\
	J_2 & := \sup_{(0,T)} \int_0^t \Vnorm{e^{(t-s) \ml_1} W_2 (u,v)}_{r_0,p_0} ds, & K_2 & := \sup_{(0,T)} t^{a_1} \int_0^t \Vnorm{e^{(t-s) \ml_1} W_2 (u,v)}_{r_1,p_1} ds, \\
	J_3 & := \sup_{(0,T)} \int_0^t \Vnorm{e^{(t-s) \ml_1} W_1 (B,B)}_{r_0,p_0} ds, & K_3 & := \sup_{(0,T)} t^{a_1} \int_0^t \Vnorm{e^{(t-s) \ml_1} W_1 (B,B)}_{r_1,p_1} ds.
\end{align*}

We will show that each term is bounded above by $C M^2 T^k$ for various values of $k > 0$, which, since $T < 1$, will imply $I_2, I_3 < M/4$ provided that $CM^2 < M/4$.

We begin our algorithm with $J_1$ and $K_1$, showing the details of the calculations and highlighting the choices of parameters. The argument for the other two pairs of integrals is very similar, and the details will be omitted.

\subsection{$J_1$ and $K_1$}

By Proposition \ref{semigroup-est}, $J_1$ is bounded by
\begin{equation*}
	J_1 \leq \sup_{(0,T)} \int_0^t (t-s)^{-(r_0 - (\gamma_3^- - 1) + n/{\pi_1} - n/p_0)/\gamma_1^-} \Vnorm{W_1 (u,v)}_{\gamma_3^- - 1,\pi_1}ds
\end{equation*}
provided that $\gamma_3^- - 1 \leq r_0$ and where $\pi_1$ is an intermediate parameter that will be specified later. Now we work towards bounding $W_1(u,v)$, and an application of Proposition \ref{inverse-semigroup-est} gives us
\begin{equation*}
	\Vnorm{W_1 (u,v)}_{\gamma_3^- - 1,\pi_1} = \Vnorm{P(1 - \ml_3)^{-1} \diverg (u \otimes v)}_{\gamma_3^- - 1,\pi_1} \leq C \Vnorm{u \otimes v}_{\pi_1}.
\end{equation*}
We chose a regularity of $\gamma_3^- - 1$ when performing the semigroup estimate in order to end up with a product in zero regularity, so that we can apply Holder's inequality: by Proposition \ref{product-est}, if 
\begin{equation} \label{pi_1}
\frac{1}{\pi_1} = \frac{1}{p'} + \frac{1}{p_1},
\end{equation}
we have
\begin{equation*}
	\Vnorm{u \otimes v}_{\pi_1} \leq C \Vnorm{u}_{p'} \Vnorm{v}_{p_1} \leq C \Vnorm{u}_{p'} \Vnorm{u}_{\gamma_3^-, p_1}.
\end{equation*}
Note that Equation (\ref{pi_1}) specifies the required value of $\pi_1$. We obtain $\Vnorm{u}_{p'} \leq \Vnorm{u}_{r_0,p_0}$ by Proposition \ref{sobolev-est} if
\begin{equation} \label{J_1-p'}
	r_0 < \frac{n}{p_0} \text{ and } \frac{1}{p'} = \frac{1}{p_0} - \frac{r_0}{n};
\end{equation}
we also get $\Vnorm{u}_{\gamma_3^-,p_1} \leq \Vnorm{u}_{r_1, p_1}$ by requiring that $r_1 \geq \gamma_3^-$. Combining the two bounds gives us
\begin{equation} \label{J_1-W_1}
	\Vnorm{W_1(u,v)}_{\gamma_3^- - 1, \pi_1} \leq C \Vnorm{u \otimes v}_{\pi_1} \leq C \Vnorm{u}_{r_0,p_0} \Vnorm{u}_{r_1, p_1}.
\end{equation}
Note that Equations (\ref{pi_1}) and (\ref{J_1-p'}) give us
\begin{equation*}
	\frac{1}{\pi_1} = \frac{1}{p_0} + \frac{1}{p_1} - \frac{r_0}{n}.
\end{equation*}

With this new bound on $W_1(u,v)$, we come back to $J_1$ and see that
\begin{align*}
	J_1 & \leq \sup_{(0,T)} \int_0^t (t-s)^{-(r_0 - (\gamma_3^- - 1) + n/{\pi_1} - n/p_0)/\gamma_1^-} \Vnorm{W_1 (u,v)}_{\gamma_3^- - 1,\pi_1}ds \\
	& \leq C \sup_{(0,T)} \int_0^t (t-s)^{-(r_0 - (\gamma_3^- - 1) + n/{\pi_1} - n/p_0)/\gamma_1^-} \Vnorm{u}_{r_0, p_0} \Vnorm{u}_{r_1, p_1} ds \\
	& = C \sup_{(0,T)} \int_0^t (t-s)^{-(r_0 - (\gamma_3^- - 1) + n/{\pi_1} - n/p_0)/\gamma_1^-} s^{-a_1} \Vnorm{u}_{r_0, p_0} s^{a_1} \Vnorm{u}_{r_1, p_1} ds \\
	& \leq C \Vnorm{u}_{0; r_0, p_0} \Vnorm{u}_{a_1; r_1, p_1} \sup_{(0,T)} \int_0^t (t-s)^{-(r_0 - (\gamma_3^- - 1) + n/{\pi_1} - n/p_0)/\gamma_1^-} s^{-a_1} ds \\
	& < CM^2 T^{1 - (r_0 - (\gamma_3^- - 1) + n/{\pi_1} - n/p_0)/\gamma_1^- - a_1},
\end{align*}
where the last inequality holds by Proposition \ref{integral-est} if
\begin{align} \label{J_1-gamma_1}
	\gamma_1^- & > r_0 - (\gamma_3^- - 1) + \frac{n}{\pi_1} - \frac{n}{p_0} + \gamma_1 a_1 \notag \\
	& = r_0 - (\gamma_3^- - 1) + n \left( \frac{1}{p_0} + \frac{1}{p_1} - \frac{r_0}{n} \right) - \frac{n}{p_0} + r_1 - r_0 + \frac{n}{p_0} - \frac{n}{p_1} \notag \\
	& = 1 - r_0 + r_1 - \gamma_3^- + \frac{n}{p_0}.
\end{align}
We further note that the requirement in (\ref{J_1-gamma_1}) also guarantees that the exponent on $T$ is positive, as desired.

We now turn our attention to $K_1$. Proposition \ref{semigroup-est} guarantees that, if $p_1 \geq \pi_1'$,
\begin{equation*}
	K_1 \leq \sup_{(0,T)} \int_0^t (t-s)^{-(r_1 - (r_1 - r_0 + \gamma_3^- - 1) + n/{\pi_1'} - n/p_1)/\gamma_1^-} \Vnorm{W_1 (u,v)}_{r_1 - r_0 + \gamma_3^- - 1,\pi_1'} ds.
\end{equation*}
This time, we chose $r_1 - r_0 + \gamma_3^- - 1$ in order to match the previous ``jump'' in regularity from $r_0$ to $\gamma_3^- - 1$. Propositions \ref{product-est} and \ref{inverse-semigroup-est} give us
\begin{align*}
	\Vnorm{W_1 (u,v)}_{r_1 - r_0 + \gamma_3^- - 1,\pi_1'} & = \Vnorm{P(1 - \ml_3)^{-1} \diverg (u \otimes v)}_{r_1 - r_0 + \gamma_3^- - 1,\pi_1'} \\
	& \leq C \Vnorm{u \otimes v}_{r_1 - r_0, \pi_1'} \\
	& \leq C \left( \Vnorm{u}_{r_1 - r_0, p'} \Vnorm{v}_{p''} + \Vnorm{v}_{r_1 - r_0, q'} \Vnorm{u}_{q''} \right), 
\end{align*}
provided that
\begin{equation*}
	\frac{1}{\pi_1'} = \frac{1}{p'} + \frac{1}{p''} = \frac{1}{q'} + \frac{1}{q''}.
\end{equation*}
Four applications of Proposition \ref{sobolev-est} lead us to the following bounds:
\begin{align}
	& \Vnorm{u}_{r_1 - r_0, p'} \leq \Vnorm{u}_{r_1,p_1} & \text{if } r_0 < \frac{n}{p_1} \text{ and } \frac{1}{p'} = \frac{1}{p_1} - \frac{r_0}{n}, \label{K_1-p'} \\
	& \Vnorm{v}_{p''} \leq \Vnorm{v}_{r_0 - \gamma_3^-, p_0} & \text{if } r_0 < \frac{n}{p_0} + \gamma_3^- \text{ and } \frac{1}{p''} = \frac{1}{p_0} - \frac{r_0 - \gamma_3^-}{n}, \label{K_1-p''} \\
	& \Vnorm{v}_{r_1 - r_0, q'} \leq \Vnorm{v}_{r_1 - \gamma_3^-, p_1} & \text{if } r_0 < \frac{n}{p_0} + \gamma_3^- \text{ and } \frac{1}{q'} = \frac{1}{p_1} - \frac{r_0 - \gamma_3^-}{n}, \label{K_1-q'} \\
	& \Vnorm{u}_{q''} \leq \Vnorm{u}_{r_0,p_0} & \text{if } r_0 < \frac{n}{p_0} \text{ and } \frac{1}{q''} = \frac{1}{p_0} - \frac{r_0}{n}. \label{K_1-q''}
\end{align}
Combining the parameters specified by (\ref{K_1-p'})-(\ref{K_1-q''}), we obtain
\begin{equation*}
	\frac{1}{\pi_1'} = \frac{1}{p_0} + \frac{1}{p_1} - \frac{2r_0 - \gamma_3^-}{n}
\end{equation*}
and
\begin{equation*}
	\Vnorm{W_1(u,v)}_{r_1 - r_0 + \gamma_3^- - 1,\pi_1'} \leq C \Vnorm{u \otimes v}_{r_1 - r_0, \pi_1'} \leq C \Vnorm{u}_{r_0,p_0} \Vnorm{u}_{r_1,p_1}.
\end{equation*}
Moreover, the integrability requirement from Proposition \ref{semigroup-est} necessitates
\begin{equation*}
	\frac{1}{p_1} \leq \frac{1}{\pi_1'} = \frac{1}{p_0} + \frac{1}{p_1} - \frac{2r_0 - \gamma_3^-}{n},
\end{equation*}
and so
\begin{equation} \label{extra-cond}
	r_0 \leq \frac{1}{2} \left( \frac{n}{p_0} + \gamma_3^- \right).
\end{equation}

We can finally plug this bound into $K_1$:
\begin{align*}
	K_1 & \leq \sup_{(0,T)} t^{a_1} \int_0^t (t-s)^{-(r_1 - (r_1 - r_0 + \gamma_3^- - 1) + n/{\pi_1'} - n/p_1)/\gamma_1^-} \Vnorm{W_1 (u,v)}_{r_1 - r_0 + \gamma_3^- - 1,\pi_1'} ds \\
	& \leq C \sup_{(0,T)} t^{a_1} \int_0^t (t-s)^{-(r_1 - (r_1 - r_0 + \gamma_3^- - 1) + n/{\pi_1'} - n/p_1)/\gamma_1^-} \Vnorm{u}_{r_0,p_0} \Vnorm{u}_{r_1,p_1} ds \\
	& = C \sup_{(0,T)} t^{a_1} \int_0^t (t-s)^{-(r_1 - (r_1 - r_0 + \gamma_3^- - 1) + n/{\pi_1'} - n/p_1)/\gamma_1^-} s^{-a_1} \Vnorm{u}_{r_0,p_0} s^{a_1} \Vnorm{u}_{r_1,p_1} ds \\
	& \leq C \Vnorm{u}_{0;r_0,p_0} \Vnorm{u}_{a_1; r_1,p_1} \sup_{(0,T)} t^{a_1}\int_0^t (t-s)^{-(r_1 - (r_1 - r_0 + \gamma_3^- - 1) + n/{\pi_1'} - n/p_1)/\gamma_1^-} s^{-a_1} ds \\
	& < CM^2 T^{1 - (r_1 - (r_1 - r_0 + \gamma_3^- - 1) + n/{\pi_1'} - n/p_0)/\gamma_1^-},
\end{align*}
where the last inequality follows by Proposition \ref{integral-est} if
\begin{align}
	\gamma_1^- & > r_1 - \left( r_1 - r_0 + \gamma_3^- - 1 \right) + \frac{n}{\pi_1'} - \frac{n}{p_1} + \gamma_1^- a_1 \notag \\
	& = r_0 - \gamma_3^-+  1 + n \left( \frac{1}{p_0} + \frac{1}{p_1} - \frac{2r_0 - \gamma_3^-}{n} \right) - \frac{n}{p_1} + r_1 - r_0 + \frac{n}{p_0} - \frac{n}{p_1} \notag \\
	& = 1 - 2r_0 + r_1 + \frac{2n}{p_0} - \frac{n}{p_1} \label{K_1-gamma_1-new}. 
\end{align}
Note that, once again, the requirement that Proposition \ref{integral-est} hold is sufficient to guarantee that the exponent on $T$ be positive.

To summarize, here is the list of inequalities needed to obtain the desired bounds on $J_1$ and $K_1$:
\begin{align*}
	& r_0 \leq \gamma_3^- & \text{assumption}, \\
	& \gamma_3^- - 1 \leq r_0 & \text{semigroup estimate for $J_1$}, \\
	& r_1 - r_0 + \gamma_3^- - 1 \leq r_1 & \text{semigroup estimate for $K_1$} \\ 
	& r_0 < \frac{n}{p_0} & (\ref{J_1-p'}), \\
	& r_0 < \frac{n}{p_1} & (\ref{K_1-p'}), \\
	& r_0 \leq \frac{1}{2} \left( \frac{n}{p_0} + \gamma_3^- \right) & (\ref{extra-cond}), \\
	& r_0 < \frac{n}{p_0} + \gamma_3^- & (\ref{J_2-p''}), \\
	& r_1 \geq \gamma_3^- & \text{bound on $\Vnorm{u}_{\gamma_3^-,p_1}$ in $J_1$}, \\
	& \gamma_1^- > 1 - r_0 + r_1 - \gamma_3^- + \frac{n}{p_0} & (\ref{J_1-gamma_1}), \\
	& \gamma_1^- > 1 - 2r_0 + r_1 + \frac{2n}{p_0} - \frac{n}{p_1} & (\ref{K_1-gamma_1-new}).
\end{align*}
After some obvious simplifications and after noting that (\ref{K_1-gamma_1-new}) implies (\ref{J_1-gamma_1}) since
\begin{equation*}
	\underbrace{\left( 1 - 2r_0 + r_1 + \frac{2n}{p_0} - \frac{n}{p_1} \right)}_{\text{RHS of (\ref{K_1-gamma_1-new})}} - \underbrace{\left( 1 - r_0 + r_1 - \gamma_3^- + \frac{n}{p_0} \right)}_{\text{RHS of (\ref{J_1-gamma_1})}} = \gamma_3^- - r_0 + \frac{n}{p_0} - \frac{n}{p_1} \geq 0,
\end{equation*}
the list reduces to
\begin{gather*}
	\gamma_3^- - 1 \leq r_0 \leq \gamma_3^- \leq r_1, \\
	r_0 < \frac{n}{p_1}, \\
	\gamma_1^- > 1 - 2r_0 + r_1 + \frac{2n}{p_0} - \frac{n}{p_1}.
\end{gather*}

\subsection{$J_2$ and $K_2$.}

We have
\begin{equation*}
	J_2 = \sup_{(0,T)} \int_0^t \Vnorm{e^{(t-s) \ml_1} W_2 (u,v)}_{r_0,p_0} ds \leq \sup_{(0,T)} \int_0^t \Vnorm{W_2 (u,v)}_{r_0, p_0} ds. 
\end{equation*}
We now work towards bounding $W_2(u,v)$. We immediately see, thanks to Proposition \ref{inverse-semigroup-est}, that
\begin{align*}
	\Vnorm{W_2(u,v)}_{r_0, p_0} & = \Vnorm{P (1 - \ml_3)^{-1} \sum_{i = 1}^n v_i \nabla u_i}_{r_0,p_0} \\
	& \leq C \Vnorm{\sum_{i = 1}^n v_i \nabla u_i}_{r_0 - \gamma_3^-,p_0} \\
	& \leq C \sum_{i = 1}^n \Vnorm{v_i \nabla u_i}_{p_0}
\end{align*}
since $r_0 \leq \gamma_3^-$. Now the product estimate is nothing more than Holder's inequality, so if
\begin{equation*}
	\frac{1}{p_0} = \frac{1}{p'} + \frac{1}{p''}
\end{equation*}
we obtain
\begin{align*}
	\Vnorm{W_2(u,v)}_{r_0, p_0} & \leq C \sum_{i = 1}^n \Vnorm{v_i}_{p'} \Vnorm{\nabla u_i}_{p''} \\
	& \leq C \sum_{i = 1}^n \Vnorm{v}_{p'} \Vnorm{\nabla u}_{p''} \\
	& \leq C \Vnorm{u}_{\gamma_3^-, p'} \Vnorm{u}_{1, p''}.
\end{align*}

By Proposition \ref{sobolev-est} we have that
\begin{equation*}
	\Vnorm{u}_{\gamma_3^-, p'} \leq C \Vnorm{u}_{\gamma_3^- + \beta, p_1} \leq C \Vnorm{u}_{r_1,p_1} \text{ and } \Vnorm{u}_{1, p''} \leq C \Vnorm{u}_{r_1,p_1},
\end{equation*}
where the first set of inequalities requires that
\begin{equation}  \label{J_2-p'}
	0 \leq \beta < \frac{n}{p_1}, \quad \frac{1}{p'} = \frac{1}{p_1} - \frac{\beta}{n}, \quad r_1 \geq \gamma_3^- + \beta,
\end{equation}
and the second inequality requires that
\begin{equation} \label{J_2-p''}
	\frac{1}{p''} = \frac{1}{p_1} - \frac{r_1 - 1}{n} \text{ and } r_1 \geq 1.
\end{equation}

We finally obtain
\begin{equation*}
	\Vnorm{W_2(u,v)}_{r_0, p_0} \leq C \Vnorm{u}_{\gamma_3^-, p'} \Vnorm{u}_{1, p''} \leq C \Vnorm{u}_{r_1, p_1}^2.
\end{equation*}
We pause here to note that, without the presence of the space $\dot{C}_{a_1;r_1,p_1}(\R^n)$ in the definition of $X_{T,M}$, we would not be able to bound this $W_2$ term.

Returning to $J_2$, we have 
\begin{align*}
	J_2 & \leq \sup_{(0,T)} \int_0^t \Vnorm{W_2 (u,v)}_{r_0, \pi_2} ds \\
	& \leq C \sup_{(0,T)} \int_0^t \Vnorm{u}_{r_1,p_1}^2 ds \\
	& = C \sup_{(0,T)} \int_0^t s^{-2a_1} s^{a_1} \Vnorm{u}_{r_1,p_1} s^{a_1} \Vnorm{u}_{r_1,p_1} ds \\
	& \leq C \Vnorm{u}_{a_1;r_1,p_1}^2 \sup_{(0,T)} \int_0^t s^{-2a_1} ds \\
	& < CM^2 T^{1 - 2a_1},
\end{align*}
provided $2a_1 > 1$, which is equivalent to
\begin{equation} \label{J_2-gamma_1-new}
	\gamma_1^- > 2r_1 - 2r_0 + \frac{2n}{p_0} - \frac{2n}{p_1}
\end{equation}
and we recall that
\begin{equation*}
	\frac{n}{p_0} = \frac{n}{p'} + \frac{n}{p''} = \frac{2n}{p_1} - \beta - r_1 + 1.
\end{equation*}
We choose $\beta$ to be exactly
\begin{equation} \label{beta}
	\beta = 1 - r_1 - \frac{n}{p_0} + \frac{2n}{p_1},
\end{equation}
and so the two requirements in (\ref{J_2-p'}) become
\begin{equation} \label{J_2-p'-new}
	r_1 \geq 1 - \frac{n}{p_0} + \frac{n}{p_1} \text{ and } 2r_1 \geq 1 + \gamma_3^- - \frac{n}{p_0} + \frac{2n}{p_1}.
\end{equation}

Turning to $K_2$, noting that we go down to $\gamma_3^-$ instead of $r_0$, we have
\begin{align*}
	K_2 & = \sup_{(0,T)} t^{a_1} \int_0^t \Vnorm{e^{(t-s) \ml_1} W_2 (u,v)}_{r_1,p_1} ds \\
	& \leq \sup_{(0,T)} t^{a_1} \int_0^t (t-s)^{-(r_1 - \gamma_3^- + n/{p_0} - n/p_1)/\gamma_1^-} \Vnorm{W_2 (u,v)}_{\gamma_3^-,p_0}ds \\
	& \leq C \sup_{(0,T)} t^{a_1} \int_0^t (t-s)^{-(r_1 - \gamma_3^- + n/{p_0} - n/p_1)/\gamma_1^-} \Vnorm{u}_{r_1,p_1} \Vnorm{u}_{r_1,p_1} ds \\
	& = C \sup_{(0,T)} t^{a_1} \int_0^t (t-s)^{-(r_1 - \gamma_3^- + n/{p_0} - n/p_1)/\gamma_1^-} s^{-2a_1} s^{a_1} \Vnorm{u}_{r_1,p_1} s^{a_1} \Vnorm{u}_{r_1,p_1} ds \\
	& \leq C \Vnorm{u}_{a_1;r_0,p_0}^2 \sup_{(0,T)} t^{a_1} \int_0^t (t-s)^{-(r_1 - \gamma_3^- + n/{p_0} - n/p_1)/\gamma_1^-} s^{-2a_1} ds \\
	& < CM^2 T^{1 - (r_1 - \gamma_3^- + n/\pi_2 - n/p_1)/\gamma_1^- - a_1},
\end{align*}
where, by Proposition \ref{integral-est}, the last inequality holds if
\begin{align} \label{K_2-gamma_1}
	\gamma_1^- & > r_1 - \gamma_3^- + \frac{n}{p_0} - \frac{n}{p_1} + 2\gamma_1^-a_1 \notag \\
	& = 3r_1 - 2r_0 - \gamma_3^- + \frac{3n}{p_0} - \frac{3n}{p_1}.
\end{align}

Here is a summary of the inequalities needed to obtain the required bounds on $J_2$ and $K_2$:
\begin{align*}
	& r_1 \geq 1 & (\ref{J_2-p''}), \\
	& r_1 \geq 1 - \frac{n}{p_0} + \frac{n}{p_1} & (\ref{J_2-p'-new}), \\
	& 2r_1 \geq 1 + \gamma_3^- - \frac{n}{p_0} + \frac{2n}{p_1} & (\ref{J_2-p'-new}), \\
	& \gamma_1^- > 2r_1 - 2r_0 + \frac{2n}{p_0} - \frac{2n}{p_1} & (\ref{J_2-gamma_1-new}), \\
	& \gamma_1^- > 3r_1 - 2r_0 - \gamma_3^- + \frac{3n}{p_0} - \frac{3n}{p_1} & (\ref{K_2-gamma_1}).
\end{align*}
By noting that
\begin{equation*}
	\underbrace{\left(3r_1 - 2r_0 - \gamma_3^- + \frac{3n}{p_0} - \frac{3n}{p_1} \right)}_{\text{RHS of (\ref{K_2-gamma_1})}} - \underbrace{\left( 2r_1 - 2r_0 + \frac{2n}{p_0} - \frac{2n}{p_1}\right)}_{\text{RHS of (\ref{J_2-gamma_1-new})}} = r_1 - \gamma_3^- + \frac{n}{p_0} - \frac{n}{p_1} \geq 0
\end{equation*}
we conclude that (\ref{K_2-gamma_1}) implies (\ref{J_2-gamma_1-new}), and so the list reduces to
\begin{gather*}
	2r_1 \geq \max \left\{2, 1 + \gamma_3^- - \frac{n}{p_0} + \frac{2n}{p_1} \right\}, \\
	\gamma_1^- > 3r_1 - 2r_0 - \gamma_3^- + \frac{3n}{p_0} - \frac{3n}{p_1}.
\end{gather*}

\subsection{$J_3$ and $K_3$.}

Provided that $r_2 - 1 + \gamma_3^- \leq r_0$ and $\frac{1}{\pi_3} \geq \frac{1}{p_0}$, Proposition \ref{semigroup-est} gives us
\begin{align*}
	J_3 & = \sup_{(0,T)} \int_0^t \Vnorm{e^{(t-s) \ml_1} W_1 (B,B)}_{r_0,p_0} ds \\
	& \leq \sup_{(0,T)} \int_0^t (t-s)^{-(r_0 - (r_2 - 1 + \gamma_3^-) + n/{\pi_3} - n/p_0)/\gamma_1^-} \Vnorm{W_1 (B,B)}_{r_2 - 1 + \gamma_3^-,\pi_3} ds.
\end{align*}
Once again, applying Propositions \ref{product-est} and \ref{inverse-semigroup-est} gets us
\begin{align*}
	\Vnorm{W_1(B,B)}_{r_2 - 1 + \gamma_3^-, \pi_3} & = \Vnorm{P (1 - \ml_3)^{-1} \diverg (B \otimes B)}_{r_2 - 1 + \gamma_3^-, \pi_3} \\
	& \leq C \Vnorm{B \otimes B}_{r_2, \pi_3} \\
	& \leq C \Vnorm{B}_{r_2, p_2} \Vnorm{B}_{p'},
\end{align*}
where the product estimate requires that $r_2 \geq 0$ and
\begin{equation*}
	\frac{1}{\pi_3} = \frac{1}{p_2} + \frac{1}{p'}.
\end{equation*}
Provided that
\begin{equation} \label{J_3-p'}
	r_2 < \frac{n}{p_2} \text{ and } \frac{1}{p'} = \frac{1}{p_2} - \frac{r_2}{n},
\end{equation}
which combines with the previous equation to give
\begin{equation*}
	\frac{1}{\pi_3} = \frac{2}{p_2} - \frac{r_2}{n},
\end{equation*}
we can bound $\Vnorm{B}_{p'}$ by $\Vnorm{B}_{r_2,p_2}$ thanks to Proposition \ref{sobolev-est}. Thus,
\begin{equation*}
	\Vnorm{W_1(B,B)}_{r_2-1+\gamma_3^-,\pi_3} \leq C \Vnorm{B}_{r_2,p_2}^2.
\end{equation*}

Moreover, Proposition (\ref{semigroup-est}) requires that
\begin{equation*}
	\frac{1}{p_0} \leq \frac{1}{\pi_3} = \frac{2}{p_2} - \frac{r_2}{n},
\end{equation*}
which can be restated as
\begin{equation} \label{r_2-extra}
	r_2 \leq \frac{2n}{p_2} - \frac{n}{p_0}.
\end{equation}

Plugging the bound for $W_1(B,B)$ back into the integral gives us
\begin{align*}
	J_3 & \leq \sup_{(0,T)} \int_0^t (t-s)^{-(r_0 - (r_2 - 1 + \gamma_3^-) + n/{\pi_3} - n/p_0)/\gamma_1^-} \Vnorm{W_1 (B,B)}_{r_2 - 1 + \gamma_3^-,\pi_3} ds \\
	& \leq C \sup_{(0,T)} \int_0^t (t-s)^{-(r_0 - (r_2 - 1 + \gamma_3^-) + n/{\pi_3} - n/p_0)/\gamma_1^-} \Vnorm{B}_{r_2,p_2}^2 ds \\
	& \leq C \Vnorm{B}_{0;r_2,p_2}^2 \sup_{(0,T)} \int_0^t (t-s)^{-(r_0 - (r_2 - 1 + \gamma_3^-) + n/{\pi_3} - n/p_0)/\gamma_1^-} ds \\
	& < C M^2 T^{1 - (r_0 - (r_2 - 1 + \gamma_3^-) + n/{\pi_3} - n/p_0)/\gamma_1^-},
\end{align*}
where once again the last inequality holds if
\begin{align}
	\gamma_1^- & > r_0 - (r_2 - 1 + \gamma_3^-) + \frac{n}{\pi_3} - \frac{n}{p_0} \notag \\
	& = r_0 - 2r_2 - \gamma_3^- + 1 + \frac{2n}{p_2} - \frac{n}{p_0}. \label{J_3-gamma_1}
\end{align}

The same bounds for $W_1(B,B)$ work in the case of $K_3$, so that
\begin{align*}
	K_3 & = \sup_{(0,T)} t^{a_1} \int_0^t \Vnorm{e^{(t-s) \ml_1} W_1 (B,B)}_{r_1,p_1} ds \\
	& \leq \sup_{(0,T)} t^{a_1} \int_0^t (t-s)^{-(r_1 - (r_2 - 1 + \gamma_3^-) + n/\pi_3 - n/p_1)/\gamma_1^-} \Vnorm{W_1 (B,B)}_{r_2 - 1 + \gamma_3^-,\pi_3} ds \\
	& \leq C \sup_{(0,T)} t^{a_1} \int_0^t (t-s)^{-(r_1 - (r_2 - 1 + \gamma_3^-) + n/\pi_3 - n/p_1)/\gamma_1^-} \Vnorm{B}_{r_2,p_2}^2 ds \\
	& \leq C \Vnorm{B}_{0;r_2,p_2}^2 \sup_{(0,T)} t^{a_1} \int_0^t (t-s)^{-(r_1 - (r_2 - 1 + \gamma_3^-) + n/\pi_3 - n/p_1)/\gamma_1^-} ds \\
	& < C M^2 T^{1 - (r_1 - (r_2 - 1 + \gamma_3^-) + n/\pi_3 - n/p_1)/\gamma_1^- + a_1},
\end{align*}
which holds provided that
\begin{align} \label{K_3-gamma_1}
	\gamma_1^- & > r_1 - (r_2 - 1 + \gamma_3^-) + \frac{n}{\pi_3} - \frac{n}{p_1} \notag \\
	& = r_1 - 2r_2 - \gamma_3^- + 1 - \frac{n}{p_1} + \frac{2n}{p_2}.
\end{align}

What follows is the list of inequalities needed to bound $J_3$ and $K_3$ as desired:
\begin{align*}
	& r_2 - 1 + \gamma_3^- \leq r_0 & \text{semigroup estimate for $J_3$}, \\
	& r_2 - 1 + \gamma_3^- \leq r_1 & \text{semigroup estimate for $K_3$}, \\
	& r_2 \geq 0 & \text{product estimate}, \\
	& r_2 < \min \left\{ \frac{n}{p_2}, \frac{2n}{p_2} - \frac{n}{p_0} \right\} & (\ref{J_3-p'})-(\ref{r_2-extra}), \\
	& \gamma_1^- > r_0- 2r_2 - \gamma_3^- + 1 - \frac{n}{p_0} + \frac{2n}{p_2} & (\ref{J_3-gamma_1}), \\
	& \gamma_1^- > r_1 - 2r_2 - \gamma_3^- + 1 + \frac{2n}{p_2} - \frac{n}{p_1} & (\ref{K_3-gamma_1}).
\end{align*}
We see that (\ref{K_3-gamma_1}) suffices for (\ref{J_3-gamma_1}) since
\begin{equation*}
	\underbrace{\left( r_1 - 2r_2 - \gamma_3^- + 1 + \frac{2n}{p_2} - \frac{n}{p_1} \right)}_{\text{RHS of (\ref{K_3-gamma_1})}} - \underbrace{\left( r_0- 2r_2 - \gamma_3^- + 1 - \frac{n}{p_0} + \frac{2n}{p_2} \right)}_{\text{RHS of (\ref{J_3-gamma_1})}} = r_1 - r_0 + \frac{n}{p_0} - \frac{n}{p_1} \geq 0,
\end{equation*}
and so the list reduces to
\begin{gather*}
	r_2 - 1 + \gamma_3^- \leq r_0, \\
	0 \leq r_2 < \min \left\{ \frac{n}{p_2}, \frac{2n}{p_2} - \frac{n}{p_0} \right\}, \\
	\gamma_1^- > r_1 - 2r_2 - \gamma_3^- + 1 + \frac{2n}{p_2} - \frac{n}{p_1}.
\end{gather*}

\subsection{Bounding $I_4$.}

Applying Minkowski's inequality to $I_4$ gives
\begin{equation*}
	I_4 \leq L_1 + L_2,
\end{equation*}
where
\begin{align*}
	L_1 & := \sup_{(0,T)} \int_0^t \Vnorm{e^{(t-s) \ml_2} W_3(u,B) ds}_{r_2,p_2} \\
	L_2 & := \sup_{(0,T)} \int_0^t \Vnorm{e^{(t-s) \ml_2} W_3(B,u) ds}_{r_2,p_2}
\end{align*}
We can immediately note that, since $W_3$ is not symmetric, $L_1 \neq L_2$, but our techniques will give the same bound for each. So, we set $L := L_1$ and proceed to bound only $L_1$. Proposition \ref{semigroup-est} gives us
\begin{equation*}
	L \leq \sup_{(0,T)} \int_0^t (t-s)^{-(r_2 - (r_2 - 1) + n/\pi_4 - n/p_2)/\gamma_2^-} \Vnorm{W_3(u,B)}_{r_2 - 1, \pi_4} ds,
\end{equation*}
provided that $\frac{1}{\pi_4} \geq \frac{1}{p_2}$.

Continuing with $W_3(u,B)$, we obtain
\begin{equation*}
	\Vnorm{W_3(u,B)}_{r_2 - 1, \pi_4} = \Vnorm{P \diverg (u \otimes B)}_{r_2 - 1, \pi_4} \leq C \Vnorm{u \otimes B}_{r_2, \pi_4};
\end{equation*}
an application of Proposition \ref{product-est} gives us
\begin{align*}
	\Vnorm{u \otimes B}_{r_2 , \pi_4} \leq C \left( \Vnorm{u}_{r_2,p'} \Vnorm{B}_{p''} + \Vnorm{B}_{r_2,p_2} \Vnorm{u}_{q''} \right)
\end{align*}
as long as
\begin{equation*}
	\frac{1}{\pi_4} = \frac{1}{p'} + \frac{1}{p''} = \frac{1}{p_2} + \frac{1}{q''}.
\end{equation*}
We want to bound $\Vnorm{u \otimes B}_{r_2, \pi_4}$ by $\Vnorm{u}_{r_0,p_0} \Vnorm{B}_{r_2,p_2}$, which requires three applications of Proposition \ref{sobolev-est}. First, we obtain $\Vnorm{u}_{r_2,p'} \leq \Vnorm{u}_{r_0,p_0}$ if
\begin{equation} \label{L-p'}
	0 \leq r_0 - r_2 < \frac{n}{p_0} \text{ and } \frac{1}{p'} - \frac{r_2}{n} = \frac{1}{p_0} - \frac{r_0}{n}.
\end{equation}
We further get $\Vnorm{u}_{q''} \leq \Vnorm{u}_{r_0,p_0}$ provided that
\begin{equation} \label{L-q''}
	r_0 < \frac{n}{p_0} \text{ and } \frac{1}{q''} = \frac{1}{p_0} - \frac{r_0}{n}.
\end{equation}
The last embedding, $\Vnorm{B}_{p''} \leq \Vnorm{B}_{r_2,p_2}$, requires
\begin{equation} \label{L-p''}
	r_2 < \frac{n}{p_2} \text{ and } \frac{1}{p''} = \frac{1}{p_2} - \frac{r_2}{n}.
\end{equation}

Combining Equations (\ref{L-p'})-(\ref{L-p''}) together gives us
\begin{equation*}
	\frac{1}{\pi_4} = \frac{1}{p_0} + \frac{1}{p_2} - \frac{r_0}{n},
\end{equation*} 
which is required to satisfy
\begin{equation} \label{extra-cond-2}
	\frac{1}{p_2} \leq \frac{1}{\pi_4} = \frac{1}{p_0} + \frac{1}{p_2} - \frac{r_0}{n} \implies r_0 \leq \frac{n}{p_0}.
\end{equation}
This is the bound we were looking for:
\begin{equation*}
	\Vnorm{W_3(u,B)}_{r_2 - 1, \pi_4} \leq C \left( \Vnorm{u}_{r_2,p'} \Vnorm{B}_{p''} + \Vnorm{B}_{r_2,q'} \Vnorm{u}_{q''} \right) \leq C \Vnorm{u}_{r_0,p_0} \Vnorm{B}_{r_2,p_2}.
\end{equation*}
We can plug the above into $L$ and obtain
\begin{align*}
	L & \leq \sup_{(0,T)} \int_0^t (t-s)^{-(1 + n/\pi_4 - n/p_2)/\gamma_2^-} \Vnorm{W_3(u,B)}_{r_2 - 1, \pi_4} ds \\
	& \leq C \sup_{(0,T)} \int_0^t (t-s)^{-(1 + n/\pi_4 - n/p_2)/\gamma_2^-} \Vnorm{u}_{r_0,p_0} \Vnorm{B}_{r_2,p_2} ds \\
	& \leq C \Vnorm{u}_{0;r_0,p_0} \Vnorm{B}_{0;r_2,p_2} \sup_{(0,T)} \int_0^t (t-s)^{-(1 + n/\pi_4 - n/p_2)/\gamma_2^-} ds \\
	& \leq C M^2 T^{1 - (1 + n/\pi_4 - n/p_2)/\gamma_2^-},
\end{align*}
which holds if
\begin{align} \label{L-gamma_2}
	\gamma_2^- & > 1 + \frac{n}{\pi_4} - \frac{n}{p_2} \notag \\
	& = 1 - r_0 + \frac{n}{p_0}.
\end{align}

The list of inequalities necessary to bound $L$ is thus
\begin{align*}
	& 0 \leq r_0 - r_2 < \frac{n}{p_0} & (\ref{L-p'}), \\
	& r_2 < \frac{n}{p_2} & (\ref{L-p''}), \\
	& r_0 \leq \frac{n}{p_0} & (\ref{extra-cond-2}), \\
	& \gamma_2^- > 1 - r_0 + \frac{n}{p_0} & (\ref{L-gamma_2}).
\end{align*}

\subsection{Wrapping up}

On one final note, we point out that since
\begin{equation*}
	\underbrace{\left( 3r_1 - 2r_0 - \gamma_3^- + \frac{3n}{p_0} - \frac{3n}{p_1} \right)}_{\text{RHS of (\ref{K_2-gamma_1})}} - \underbrace{\left( 1 - 2r_0 + r_1 + \frac{2n}{p_0} - \frac{n}{p_1} \right)}_{\text{RHS of (\ref{K_1-gamma_1-new})}} = 2r_1 - 1 - \gamma_3^- + \frac{n}{p_0} - \frac{2n}{p_1} \geq 0
\end{equation*}
we have that (\ref{K_2-gamma_1}) implies (\ref{K_1-gamma_1-new}), and so the following is the definitive list containing all the inequalities needed for $I_2, I_3,$ and $I_4$:
\begin{gather*}
	\gamma_3^- - 1 \leq r_0 \leq \gamma_3^- \leq r_1, \\
	r_2 - 1 + \gamma_3^- \leq r_0, \\
	r_2 \leq r_0 < \frac{n}{p_1}, \\
	2r_1 \geq \max \left\{2, 1 + \gamma_3^- - \frac{n}{p_0} + \frac{2n}{p_1} \right\}, \\
	r_2 < \min \left\{\frac{n}{p_2}, \frac{2n}{p_2} - \frac{n}{p_0} \right\}, \\
	\gamma_1^- > 3r_1 - 2r_0 - \gamma_3^- + \frac{3n}{p_0} - \frac{3n}{p_1}, \\
	\gamma_1^- > 1 - 2r_2 + r_1  - \gamma_3^- - \frac{n}{p_1} + \frac{2n}{p_2}, \\
	\gamma_2^- > 1 - r_0 + \frac{n}{p_0}.
\end{gather*}
The above inequalities coincide with those in Theorem \ref{main-thm}, and so we are done.
\end{proof}

\end{document}